\def\titlep{Biideals and 
a lattice of C$^{*}$-bialgebras
associated with prime numbers}
\documentclass[11pt]{article}
\usepackage{graphicx,ifthen}
\usepackage{amssymb}

\font\germ=eufm10 at12pt

\def\goth#1{\hbox{\germ#1}}

\newcommand{\qed}{\hbox{\rule[-2pt]{3pt}{6pt}}}
\newcommand{\qedh}{\hfill\qed \\}

\newcommand{\vep}{\varepsilon}

\def\labelenumi{\theenumi}
\def\theenumi{\arabic{enumi}}
\def\labelenumi{\theenumi}
\def\theenumi{\Alph{enumi}}
\renewcommand{\theenumi}{\Alph{enumi}}
\def\labelenumi{\theenumi}
\def\theenumi{\arabic{enumi}}
\def\labelenumi{\theenumi}
\def\theenumi{{\rm (\roman{enumi})}}
\newtheorem{Thm}{Theorem}[section]

\newtheorem{rem}[Thm]{Remark}

\newtheorem{defi}[Thm]{Definition}
\newtheorem{lem}[Thm]{Lemma}

\newtheorem{prop}[Thm]{Proposition}

\newtheorem{cor}[Thm]{Corollary}
\newtheorem{fact}[Thm]{Fact}

\def\cal#1{\mathcal #1}
\def\con{{\cal O}_{n}}

\def\edot{=1,\ldots,n}
\def\pr{{\it Proof.}\quad}

\def\co#1{{\cal O}_{#1}}
\def\disp#1{{\displaystyle #1}}
\def\brl{branching law}
\def\bfsnl{{\rm BFS}_{N}(\Lambda)}

\addtocounter{footnote}{1}
\def\sftt#1{
\setcounter{equation}{0}
\addtocounter{footnote}{1}
\section{#1}
}

\def\ssft#1{\subsection{#1}}

\begin{document}
%
% Personal data
%
\def\autherp{Katsunori Kawamura}
\def\emailp{e-mail: kawamura@kurims.kyoto-u.ac.jp.}
\def\addressp{{\small {\it College of Science and Engineering Ritsumeikan University,}}\\
{\small {\it 1-1-1 Noji Higashi, Kusatsu, Shiga 525-8577, Japan}}
}

\def\infw{\Lambda^{\frac{\infty}{2}}V}
\def\zhalfs{{\bf Z}+\frac{1}{2}}
\def\ems{\emptyset}
\def\pmvac{|{\rm vac}\!\!>\!\! _{\pm}}
\def\vac{|{\rm vac}\rangle _{+}}
\def\dvac{|{\rm vac}\rangle _{-}}
\def\ovac{|0\rangle}
\def\tovac{|\tilde{0}\rangle}
\def\expt#1{\langle #1\rangle}
\def\zph{{\bf Z}_{+/2}}
\def\zmh{{\bf Z}_{-/2}}
\def\brl{branching law}
\def\bfsnl{{\rm BFS}_{N}(\Lambda)}
\def\scm#1{S({\bf C}^{N})^{\otimes #1}}
\def\mqb{\{(M_{i},q_{i},B_{i})\}_{i=1}^{N}}
\def\zhalf{\mbox{${\bf Z}+\frac{1}{2}$}}
\def\zmha{\mbox{${\bf Z}_{\leq 0}-\frac{1}{2}$}}
\newcommand{\mline}{\noindent
\thicklines
\setlength{\unitlength}{.1mm}
\begin{picture}(1000,5)
\put(0,0){\line(1,0){1250}}
\end{picture}
\par
 }
\def\ptimes{\otimes_{\varphi}}
\def\delp{\Delta_{\varphi}}
\def\delps{\Delta_{\varphi^{*}}}
\def\gamp{\Gamma_{\varphi}}
\def\gamps{\Gamma_{\varphi^{*}}}
\def\sem{{\sf M}}
\def\sen{{\sf N}}
\def\hdelp{\hat{\Delta}_{\varphi}}
\def\tilco#1{\tilde{\co{#1}}}
\def\ndm#1{{\bf M}_{#1}(\{0,1\})}
\def\fs{{\cal F}{\cal S}({\bf N})}

%
%%%%%%%%% Cut from here %%%%%%%%%%
%\input comm.txt
%%%%%%%%% End of Cut %%%%%%%%%
%
%
\setcounter{section}{0}
\setcounter{footnote}{0}
\setcounter{page}{1}
\pagestyle{plain}

%
% Title
%
%%%%%%%%%%%%%%%%%%%%%%%%%%%%%%%%%%%%%%%%%%
\title{\titlep}
\author{\autherp\thanks{\emailp}
\\
\addressp}
\date{}
\maketitle
%
% Abstract
%
\begin{abstract}
Let ${\cal O}_{*}$ be the C$^{*}$-algebra
defined as the direct sum of all Cuntz algebras.
Then ${\cal O}_{*}$ has a non-cocommutative comultiplication $\Delta_{\varphi}$ and a counit 
$\varepsilon$.
Let ${\rm BI}({\cal O}_{*})$ denote the set of all closed biideals of 
the C$^{*}$-bialgebra $({\cal O}_{*},\Delta_{\varphi},\varepsilon)$
and let ${\cal P}({\bf P})$ denote the power set of the set of all prime numbers.
We show a one-to-one correspondence between ${\rm BI}({\cal O}_{*})$ and ${\cal P}({\bf P})$.
Furthermore,
we show that for any ${\cal I}$ in ${\rm BI}({\cal O}_{*})$,
there exists a C$^{*}$-subbialgebra ${\cal B}_{{\cal I}}$ of ${\cal O}_{*}$
such that ${\cal O}_{*}={\cal B}_{{\cal I}}\oplus {\cal I}$, and 
the set of all such C$^{*}$-subbialgebras is a lattice 
with respect to the natural operations among C$^{*}$-subbialgebras,
which is isomorphic to the lattice ${\cal P}({\bf P})$.
\end{abstract}


%%%%%%%%%%%%%%%%%%%%%%%%%%%%%%%%%%%%%%%%%%%%%%%%%%%%%%%%%%
%
% Section 1
%
\sftt{Introduction}
\label{section:first}
We have studied C$^{*}$-bialgebras.
In this paper, we completely classify closed biideals of a certain C$^{*}$-bialgebra
by using the power set of the set of all prime numbers.
In this section, we show our motivation,
definitions of C$^{*}$-bialgebras and main theorems.

%%%%%%%%%%%%%%%%%%%%%%%%%%%%%%%%%%%%%%%%%%%%%%%%%%%%%%%%%%%%%%%%%%%%%%%%
%
% subsection 1.1
%
\ssft{Motivation}
\label{subsection:firstone}
In this subsection, we roughly explain our motivation and the background of this study.
Explicit mathematical definitions will be shown after $\S$ \ref{subsection:firsttwo}.
In \cite{TS02}, we constructed the C$^{*}$-bialgebra $\co{*}$
defined as the direct sum of all Cuntz algebras except $\co{\infty}$:
%
% Equation 1.1
%
\begin{equation}
\label{eqn:cuntztwo}
\co{*}=\co{1}\oplus\co{2}\oplus\co{3}\oplus\co{4}\oplus\cdots,
\end{equation}
where $\co{1}$ denotes the $1$-dimensional C$^{*}$-algebra ${\bf C}$
for convenience.
The C$^{*}$-bialgebra $\co{*}$ is non-commutative and non-cocommutative.
We investigated a Haar state, KMS states, C$^{*}$-bialgebra automorphisms,
C$^{*}$-subbialgebras and a comodule-C$^{*}$-algebra of $\co{*}$.
This study was motivated by a certain tensor product 
of representations of Cuntz algebras \cite{TS01}.
With respect to the tensor product, 
tensor products of irreducible representations 
and those of type III factor representations were computed \cite{TS01,TS07}.
Since there is no standard comultiplication of Cuntz algebras,
$\co{*}$ is not a deformation of any known cocommutative bialgebra. 
The C$^{*}$-bialgebra $\co{*}$ is 
a rare example of not only C$^{*}$-bialgebra but also purely algebraic bialgebra.
In fact, there exists a dense $*$-subbialgebra ${\cal A}_{0}$ of $\co{*}$
such that the image of ${\cal A}_{0}$ by the comultiplication
is contained in the algebraic tensor square of ${\cal A}_{0}$.
Furthermore, the construction of ${\cal A}_{0}$ is unlike any known
\cite{EAbe,Kassel,Sweedler}.
Hence we are interested in the bialgebra structure of $\co{*}$. 

On the other hand, 
the study of ideals of an algebra $A$ is an important basis 
for the study of $A$ itself \cite{Tomforde}.
In a similar fashion,
the study of biideals of a bialgebra $B$ will be also important 
for the study of $B$ itself.
Unfortunately, studies of biideals are few \cite{EAbe,FA,SS,Sweedler}.
Hence it is worth to construct nontrivial examples of  biideal
of a given bialgebra
and classify them. 

In this paper,
we closely consider biideals of the C$^{*}$-bialgebra $\co{*}$.

%%%%%%%%%%%%%%%%%%%%%%%%%%%%%%%%%%%%%%%%%%%%%%%%%%%%%%%%%%%%%%%%%%%%%%%%
%
% subsection 1.2
%
\ssft{C$^{*}$-bialgebra $(\co{*},\delp,\vep)$}
\label{subsection:firsttwo}
In this subsection, we recall the C$^{*}$-bialgebra in \cite{TS02}.
At first,
we prepare terminologies about C$^{*}$-bialgebra.
For two C$^{*}$-algebras $A$ and $B$,
we write ${\rm Hom}(A,B)$ the set of all $*$-homomorphisms from $A$ to $B$.
Assume that every tensor product $\otimes$ as below means 
the minimal C$^{*}$-tensor product.
A pair $(A,\Delta)$ is a {\it C$^{*}$-bialgebra}
if $A$ is a C$^{*}$-algebra and $\Delta\in {\rm Hom}(A,M(A\otimes A))$ 
where $M(A\otimes A)$ denotes the multiplier algebra of $A\otimes A$
such that the linear span of $\{\Delta(a)(b\otimes c):a,b,c\in A\}$ 
is norm dense in $A\otimes A$, and the following holds:
%
% Equation 1.2
%
\begin{equation}
\label{eqn:bialgebratwo}
(\Delta\otimes id)\circ \Delta=(id\otimes\Delta)\circ \Delta.
\end{equation}
We call $\Delta$ the {\it comultiplication} of $A$.
About C$^{*}$-bialgebras in quantum groups, see \cite{KV,MNW}.
A C$^{*}$-bialgebra $(A,\Delta)$ is {\it counital}
if there exists $\vep\in {\rm Hom}(A,{\bf C})$ such that
%
% Equation 1.3
%
\begin{equation}
\label{eqn:counit}
(\vep\otimes id)\circ \Delta= id = (id\otimes \vep)\circ
\Delta.
\end{equation}
We call $\vep$ the {\it counit} of $A$ and write
$(A,\Delta,\vep)$ as the counital C$^{*}$-bialgebra $(A,\Delta)$
with the counit $\vep$.
We state that a C$^{*}$-bialgebra $(A,\Delta)$ is {\it strictly proper}
if $\Delta(a)\in A\otimes A$ for any $a\in A$.
For two strictly proper counital 
C$^{*}$-bialgebras $(A_{1},\Delta_{1},\vep_{1})$ and 
$(A_{2},\Delta_{2},\vep_{2})$,
$f$ is a {\it strictly proper counital C$^{*}$-bialgebra morphism} from 
$(A_{1},\Delta_{1},\vep_{1})$ to $(A_{2},\Delta_{2},\vep_{2})$ 
if $f$ is $*$-homomorphism from $A_{1}$ to $A_{2}$
such that $(f\otimes f)\circ \Delta_{1}=\Delta_{2}\circ f$
and $\vep_{2}\circ f=\vep_{1}$.
A strictly proper counital C$^{*}$-bialgebra $(A_{0},\Delta_{0},\vep_{0})$ 
is a {\it strictly proper counital C$^{*}$-subbialgebra} of 
a counital strictly proper C$^{*}$-bialgebra $(A,\Delta,\vep)$ 
if $A_{0}$ is a C$^{*}$-subalgebra of $A$ such that 
$\Delta_{0}=\Delta|_{A_{0}}$ and $\vep_{0}=\vep|_{A_{0}})$.

A subspace ${\cal I}$ of a counital C$^{*}$-bialgebra $(A,\Delta,\vep)$
is a {\it biideal} of $(A,\Delta,\vep)$ \cite{EAbe,Sweedler}
if 
${\cal I}$ is a two-sided ideal of $A$ and 
%
% Equation 1.4
%
\begin{equation}
\label{eqn:biideal}
\Delta({\cal I})\subset {\cal I}\otimes A+A\otimes {\cal I},\quad
\vep({\cal I})=0. 
\end{equation}
In addition, if ${\cal I}$ is a closed subspace of $A$,
then we call ${\cal I}$ a {\it closed biideal}. 
From now,
the equivalence of subialgebras or biideals of C$^{*}$-bialgebra
is defined by strictly proper counital C$^{*}$-bialgebra morphisms.
In this paper, we treat only 
strictly proper counital C$^{*}$-bialgebras and
strictly proper counital C$^{*}$-subbialgebras.
%
% Remark 1.1
%
\begin{rem}
\label{rem:def}
{\rm
\begin{enumerate}
%(i)
\item
The notion of biideal is defined as a two-sided ideal 
and coideal \cite{EAbe,Sweedler}
even if there is no adjective ``two-sided".
%(ii)
\item
If ${\cal I}$ is a closed two-sided ideal of a C$^{*}$-algebra,
then $\{x^{*}:x\in {\cal I}\}={\cal I}$ (\cite{Ped}, Corollary 1.5.3).
Hence a closed biideal of a C$^{*}$-bialgebra
is also closed with respect to the $*$-operation.
%(iii)
\item
For a strictly proper counital C$^{*}$-bialgebra $(A,\Delta,\vep)$,
$\vep$ is not the zero map on $A$ from (\ref{eqn:counit}).
Hence $A$ itself is not a biideal of $A$ because 
of the second equation in (\ref{eqn:biideal}).
Especially, $\{0\}$ is a unique ``trivial" biideal.
%(iv)
\item
For a strictly proper counital C$^{*}$-bialgebra $(A,\Delta,\vep)$,
if ${\cal I}$ is a closed biideal of $A$,
then the quotient C$^{*}$-algebra $A/{\cal I}$ is also
a counital C$^{*}$-bialgebra 
with respect to the pair $(\tilde{\Delta},\tilde{\vep})$ defined by
$\tilde{\Delta}(x+{\cal I})\equiv \Delta(x)+{\cal I}$
and 
$\tilde{\vep}(x+{\cal I})\equiv \vep(x)+{\cal I}$
for $x\in A$.
Especially, the natural projection from $A$ to $A/{\cal I}$
is a $*$-bialgebra morphism.
%(v)
\item
The notion ``strictly proper" was introduced in \cite{TS02}
in order to treat $\co{*}$ like the purely algebraic theory of bialgebras 
\cite{EAbe,Kassel}.
This condition is quite natural from a stand point of the purely algebraic theory,
but a C$^{*}$-bialgebra is not strictly proper 
and it does not have a counit in general, 
because of the motivation of the study (\cite{KV}, p.547, ``Quantizing Locally Compact Groups").
\end{enumerate}
}
\end{rem}

Next, we introduce a C$^{*}$-bialgebra. 
Let $\con$ denote the Cuntz algebra for $2\leq n<\infty$ \cite{C},
that is, the C$^{*}$-algebra which is universally generated by
generators $s_{1},\ldots,s_{n}$ satisfying
$s_{i}^{*}s_{j}=\delta_{ij}I$ for $i,j\edot$ and
$\sum_{i=1}^{n}s_{i}s_{i}^{*}=I$
where $I$ denotes the unit of $\con$.
The Cuntz algebra $\con$ is simple, that is,
there is no nontrivial closed two-sided ideal.

Define the C$^{*}$-algebra $\co{*}$
as the direct sum of the family $\{\con:n\in {\bf N}\}$ of Cuntz algebras:
%
% Equation 1.5
%
\begin{equation}
\label{eqn:cuntbi}
\co{*}\equiv \bigoplus_{n\in {\bf N}} \con
=\{(x_{n}):\|(x_{n})\|\to 0\mbox{ as }n\to\infty\}
\end{equation}
where ${\bf N}=\{1,2,3,\ldots\}$
and $\co{1}$ denotes the $1$-dimensional C$^{*}$-algebra for convenience.
For $n\in {\bf N}$,
let $I_{n}$ denote the unit of $\con$ 
and let $s_{1}^{(n)},\ldots,s_{n}^{(n)}$ denote
canonical generators of $\con$
where $s_{1}^{(1)}\equiv I_{1}$.
For $n,m\in {\bf N}$,
define the embedding $\varphi_{n,m}$ of $\co{nm}$
into $\con\otimes \co{m}$ by
%
% Equation 1.6
%
\begin{equation}
\label{eqn:embeddingone}
\varphi_{n,m}(s_{m(i-1)+j}^{(nm)})\equiv s_{i}^{(n)}\otimes s_{j}^{(m)}
\quad(i=1,\ldots,n,\,j=1,\ldots,m).
\end{equation}
%
% Theorem 1.2
%
\begin{Thm}(\cite{TS02}, Theorem 1.1)
\label{Thm:mainone}
For the set $\varphi\equiv \{\varphi_{n,m}:n,m\in {\bf N}\}$ in
(\ref{eqn:embeddingone}),
define  $\delp^{(n)}\in {\rm Hom}(\con,\co{*}\otimes \co{*})$,
$\delp\in {\rm Hom}(\co{*},\co{*}\otimes \co{*})$ 
and $\vep\in {\rm Hom}(\co{*},{\bf C})$ by
%
% Equation 1.7,1.8,1.9
%
\begin{eqnarray}
\label{eqn:dpone}
\delp\equiv& \oplus\{\delp^{(n)}:n\in {\bf N}\},\\
\nonumber
\\
\label{eqn:dptwo}
\delp^{(n)}(x)\equiv &\disp{\sum_{(m,l)\in {\bf N}^{2},\,ml=n}\varphi_{m,l}(x)}
\quad(x\in \con,\,n\in {\bf N}),\\
\nonumber
\\
\label{eqn:dpthree}
\vep(x)\equiv  &
\left\{
\begin{array}{ll}
0\quad &\mbox{ when }x\in \oplus \{\con:n\geq 2\},\\
\\
x\quad&\mbox{ when }x\in \co{1}.
\end{array}
\right.
\end{eqnarray}
Then $(\co{*},\delp,\vep)$ is 
a strictly proper non-cocommutative counital C$^{*}$-bialgebra.
\end{Thm}
About properties of $\co{*}$, see \cite{TS02}.
About a generalization of $\co{*}$, see \cite{TS05}.

%%%%%%%%%%%%%%%%%%%%%%%%%%%%%%%%%%%%%%%%%%%%%%%%%%%%%%%%%
%
% Subsection 1.2
%
\ssft{C$^{*}$-subbialgebras and closed biideals}
\label{subsection:firstthree}
In this subsection, we consider 
C$^{*}$-subbialgebras and closed biideals 
of $(\co{*},\delp,\vep)$ in Theorem \ref{Thm:mainone}.

Let ${\cal P}({\bf N})$ denote the power set of ${\bf N}$.
For $S\in {\cal P}({\bf N})$, define the closed two-sided ideal 
$\co{*}(S)$ of the C$^{*}$-algebra $\co{*}$ by
%
% Equation 1.10
%
\begin{equation}
\label{eqn:ostar}
\co{*}(S)\equiv \oplus\{\con:n\in S\}\quad(S\ne\emptyset),\quad
\co{*}(\emptyset)\equiv 0.
\end{equation}
%
%
% Remark 1.3
%
\begin{rem}
\label{rem:not}
{\rm
From the simplicity of $\con$ and the definition of $\co{*}$,
the set of all closed two-sided ideals of $\co{*}$ coincides 
with $\{\co{*}(S):S\in {\cal P}({\bf N})\}$.
Since $\con\cong \co{m}$ if and only if $n=m$,
$\co{*}(S)\cong \co{*}(T)$ if and only if $S=T$.
In this way, the classification of closed two-sided ideals of 
the C$^{*}$-algebra $\co{*}$ is very clear.
{\it However},
for $S\in {\cal P}({\bf N})$, 
the C$^{*}$-subalgebra $\co{*}(S)$ in (\ref{eqn:ostar}) is 
{\it neither} a C$^{*}$-subbialgebra {\it nor} a biideal 
of $(\co{*},\delp,\vep)$ in general.
For example,
let 
$S\equiv \{4^{n}:n\geq 0\}$.
Then we see that $\delp(\co{*}(S))\not \subset \co{*}(S)\otimes \co{*}(S)$
by the definition of $\delp$ (see (\ref{eqn:htwo})).
Therefore 
neither
the classification of C$^{*}$-subbialgebras nor that of
biideals of $\co{*}$ is trivial. 
}
\end{rem}

A {\it monoid} is a semigroup with a unit.
A {\it submonoid} is a nonempty subset of a monoid $\sem$
which is closed with respect to the operation of $\sem$
and contains the unit of $\sem$.
We regard ${\bf N}$ as an abelian monoid with respect to the multiplication.
We write this as $({\bf N},\cdot)$.
The first main theorem in this paper is given as follows.
%
% Theorem 1.4
%
\begin{Thm}
\label{Thm:maintwo}
Let $(\co{*},\delp,\vep)$ be as in Theorem \ref{Thm:mainone}.
For any submonoid $H$ of $({\bf N},\cdot)$,
define
%
% Equation 1.11, 1.12
%
\begin{eqnarray}
\label{eqn:delpone}
\Delta_{\varphi^{(H)}}\equiv &\oplus\{\Delta_{\varphi^{(H)}}^{(n)}:n\in H\},\\
\nonumber	
\\
\label{eqn:delptwo}
\disp{\Delta_{\varphi^{(H)}}^{(n)}(x)\equiv }
&\disp{\sum_{(m,l)\in H^{2},\,ml=n}\varphi_{m,l}(x)
\quad(x\in \con,\,n\in H).}
\end{eqnarray}
Then  $(\co{*}(H),\Delta_{\varphi^{(H)}},\vep|_{{\co{*}(H)}})$
is a strictly proper counital C$^{*}$-bialgebra.
\end{Thm}

\noindent
Remark that $\Delta_{\varphi^{(H)}}\ne \delp|_{\co{*}(H)}$ in general.
Hence $(\co{*}(H),\Delta_{\varphi^{(H)}},\vep|_{{\co{*}(H)}})$
is not a C$^{*}$-subbialgebra of $(\co{*},\delp,\vep)$
in general.
For example,
$H=\{4^{l}:l\geq 0\}$ is a submonoid of $({\bf N},\cdot)$ and 
%
% Equation 1.13, 1.14
%
\begin{eqnarray}
\label{eqn:hone}
\Delta_{\varphi^{(H)}}(s_{1}^{(4)})=&I_{1}\otimes s_{1}^{(4)}+
s_{1}^{(4)}\otimes I_{1},\\
\nonumber
\\
\label{eqn:htwo}
\delp(s_{1}^{(4)})=&I_{1}\otimes s_{1}^{(4)}+
s_{1}^{(2)}\otimes s_{1}^{(2)}+s_{1}^{(4)}\otimes I_{1}.
\end{eqnarray}
Since $\Delta_{\varphi^{(H)}}\ne \delp|_{\co{*}(H)}$ in this case,
$(\co{*}(H),\Delta_{\varphi^{(H)}},\vep|_{{\co{*}(H)}})$
is not a C$^{*}$-subbialgebra of $(\co{*},\delp,\vep)$.

We consider a sufficient condition
that $\co{*}(S)$ is a C$^{*}$-subbialgebra of $(\co{*},\delp,\vep)$.
For this purpose,
we introduce several notions of monoid according to \cite{Howie,Lothaire}.
For a nonempty subset $S$ of a monoid $\sem$,
we state that $S$ is {\it proper} if $S\ne \sem$;
$S$ is a {\it subsemigroup} if $ab\in S$ when $a,b\in S$;
$S$ is {\it factorial} if $S\ne \sem$ and $a,b\in S$
when $a,b\in \sem$ and $ab\in S$;
$S$ is {\it prime} if $S\ne \sem$ and $a\in S$ or $b\in S$
when $a,b\in \sem$ and $ab\in S$;
$S$ is an {\it ideal} if $aS$, $Sa\subset  S$ for any $a\in \sem$
where $aS\equiv \{as:s\in S\}$ and $Sa\equiv \{sa:s\in S\}$.
Then we can state the next main theorem as follows.
%
% Theorem 1.5
%
\begin{Thm}
\label{Thm:mainfours}
Let $(\co{*},\delp,\vep)$ be as in Theorem \ref{Thm:mainone}.
\begin{enumerate}
%(i)
\item
For a subspace ${\cal I}$ of $\co{*}$,
${\cal I}$ is a closed biideal 
of $(\co{*},\delp,\vep)$ if and only if
\begin{enumerate}
%(a)
\item
${\cal I}=\{0\}$, or 
%(b)
\item
there exists a prime ideal ${\goth a}$ of $({\bf N},\cdot)$
such that ${\cal I}=\co{*}({\goth a})$.
\end{enumerate}
%(iii)
\item
For any closed biideal ${\cal I}$ of $(\co{*},\delp,\vep)$,
there exists a C$^{*}$-subbialgebra ${\cal B}_{{\cal I}}$ 
of $(\co{*},\delp,\vep)$
such that the quotient C$^{*}$-bialgebra 
$\co{*}/{\cal I}$ is isomorphic to ${\cal B}_{{\cal I}}$
and the following decomposition holds:
%
% Equation 1.15
%
\begin{equation}
\label{eqn:decomposition}
\co{*}={\cal B}_{{\cal I}}\oplus {\cal I}
\end{equation}
where $\oplus$ means the direct sum  of two C$^{*}$-subalgebras.
\end{enumerate}
\end{Thm}

\noindent
From Theorem \ref{Thm:mainfours}(i) and (\ref{eqn:ostar}),
for two closed biideals ${\cal I}_{1}$ and  ${\cal I}_{2}$ of $\co{*}$,
if ${\cal I}_{1}\ne {\cal I}_{2}$,
then ${\cal I}_{1}$ and ${\cal I}_{2}$ are not isomorphic as closed
two-sided ideals of $\co{*}$.
Hence ${\cal I}_{1}$ and ${\cal I}_{2}$ are not isomorphic as closed
biideals of $\co{*}$ if  ${\cal I}_{1}\ne {\cal I}_{2}$.

From Theorem \ref{Thm:mainfours}(i) and (\ref{eqn:ostar}), the following holds.
%
% Corollary 1.6
%
\begin{cor}
\label{cor:order}
Let ${\rm Spec}{\bf N}$ denote the set of all prime ideals of $({\bf N},\cdot)$.
Define $\overline{{\rm Spec}{\bf N}}\equiv {\rm Spec}{\bf N}\cup \{\emptyset\}$
and let ${\rm BI}(\co{*})$ denote the set of all closed biideals of $\co{*}$.
Then the following isomorphism between two ordered sets holds:
%
% Equation 1.16
%
\begin{equation}
\label{eqn:orderone}
\overline{{\rm Spec}{\bf N}}\ni x\longmapsto \co{*}(x)\in {\rm BI}(\co{*})
\end{equation}
where the order of both 
$\overline{{\rm Spec}{\bf N}}$ and ${\rm BI}(\co{*})$
is taken as the inclusion of subsets.
\end{cor}
Remark that the unique maximal prime ideal of ${\bf N}$
is ${\bf N}\setminus\{1\}$ and 
any minimal one is $p{\bf N}$ for a prime number $p$.
The minimal element in $\overline{{\rm Spec}{\bf N}}$ is ``$\emptyset$."
On the other hand,
the unique minimal element in ${\rm BI}(\co{*})$ is $\{0\}$
and the unique maximal one is $\bigoplus_{n\geq 2}\con$.

In order to classify  closed biideals of $\co{*}$ in $\S$ \ref{subsection:firstfour},
we prepare the following fact.
%
% Fact 1.6
%
\begin{fact}
\label{fact:inverse}
Let $\sem$ be a monoid.
For a subset $S$ of $\sem$,
$S$ is a factorial submonoid if and only if
$S^{c}\equiv \{x\in\sem:x\not\in S\}$ is a prime ideal.
\end{fact}

%%%%%%%%%%%%%%%%%%%%%%%%%%%%%%%%%%%%%%%%%%%%%%%%%%%%%%%%
%
% Subsection 1.3
%
\ssft{Lattice of C$^{*}$-subbialgebras}
\label{subsection:firstfour}
In this subsection,
we show the biideal structure of $\co{*}$ at great length.
From Theorem \ref{Thm:mainfours}(i) and Fact \ref{fact:inverse}, 
the classification of all nonzero closed biideals of $(\co{*},\delp,\vep)$
is equivalent to that of
all factorial submonoids of $({\bf N},\cdot)$.
In this subsection, we classify factorial submonoids of ${\bf N}$
and consider a correspondence between
nonzero closed biideals of $\co{*}$ and prime ideals of ${\bf N}$.

Let ${\bf P},{\cal P}({\bf P})$ and $\fs$
denote the set of all prime numbers,
the power set of ${\bf P}$ and 
the set of all factorial submonoids of ${\bf N}$, respectively.
For $F\in {\cal P}({\bf P})$,
if $F\ne \emptyset$,
define $[F]$ the submonoid of $({\bf N},\cdot)$ generated by $F$,
that is, the smallest submonoid containing $F$,
 and 
define $[\emptyset]\equiv \{1\}$.
Since ${\bf N}$ is isomorphic to the free abelian 
monoid generated by ${\bf P}$,
$H$ is a factorial submonoid of ${\bf N}$ if and only if 
there  exists $F\in {\cal P}({\bf P})\setminus\{{\bf P}\}$ such that $H=[F]$.
Define two operations $\vee$ and $\wedge$ on $\overline{\fs}\equiv \fs\cup\{{\bf N}\}$ by
%
% Equation 1.17
%
\begin{equation}
\label{eqn:operations}
H\vee G\equiv [ H\cup G],\quad 
H\wedge G\equiv H\cap G\quad(G,H\in \overline{\fs}).
\end{equation}
Then the following one-to-one correspondence
between the lattice 
$({\cal P}({\bf P}),\cup,\cap)$ 
and $(\overline{\fs} ,\vee,\wedge)$ is a lattice isomorphism \cite{Gratzer,Halmos}:
%
% Equation 1.18
%
\begin{equation}
\label{eqn:powers}
{\cal P}({\bf P})\ni F\mapsto [F]\in \overline{\fs}.
\end{equation}

For $F\in {\cal P}({\bf P})\setminus\{{\bf P}\}$, $[F]^{c}\in {\rm Spec}{\bf N}$
by Fact \ref{fact:inverse}.
For $F\in {\cal P}({\bf P})$, 
define the C$^{*}$-subbialgebra ${\cal A}(F)$ 
and the closed biideal ${\cal I}(F)$ of $\co{*}$ by
%
% Equation 1.19
%
\begin{equation}
\label{eqn:af}
{\cal A}(F)\equiv \co{*}([F]),\quad{\cal I}(F)\equiv \co{*}([F]^{c}).
\end{equation}
By definition,
the following holds for $F,G\in {\cal P}({\bf P})$:
%
% Equation 1.20
%
\begin{equation}
\label{eqn:inc}
{\cal A}(F)\subset  {\cal A}(G)\mbox{ when } F\subset G,
\end{equation}
%
% Equation 1.21
%
\begin{equation}
\label{eqn:equiv}
{\cal A}(F)={\cal A}(G)\mbox{ if and only if } F= G.
\end{equation}
Then (\ref{eqn:decomposition}) is concretely written as follows: 
%
% Equation 1.22
%
\begin{equation}
\label{eqn:exact}
\co{*}={\cal A}(F)\oplus {\cal I}(F)
\quad(F\in {\cal P}({\bf P})).
\end{equation}
Define operations $\vee,\wedge$ 
on the family ${\cal F}\equiv \{{\cal A}(F):F\in {\cal P}({\bf P})\}$ by
%
% Equation 1.23
%
\begin{equation}
\label{eqn:operationstwo}
{\cal B}\vee {\cal C}\equiv C^{*}\langle {\cal B}\cup {\cal C}\rangle,\quad
{\cal B}\wedge {\cal C}\equiv {\cal B}\cap {\cal C}
\end{equation}
for ${\cal B},{\cal C}\in {\cal F}$
where 
$C^{*}\langle X\rangle$ denotes the C$^{*}$-subalgebra of $\co{*}$ generated 
by a subset $X$ of $\co{*}$.
Then $({\cal F},\vee,\wedge)$ is a lattice.
%
% Theorem 1.8
%
\begin{Thm}
\label{Thm:mainfive}
\begin{enumerate}
%(i)
\item
For a subspace ${\cal J}$ of $\co{*}$,
${\cal J}$ is a closed biideal 
of $(\co{*},\delp,\vep)$ if and only if
there exists $F\in {\cal P}({\bf P})$ such that
${\cal J}={\cal I}(F)$.
%(ii)
\item
The lattice $\{{\cal A}(F):F\in {\cal P}({\bf P})\}$
is isomorphic to the lattice ${\cal P}({\bf P})$.
Furthermore,
${\cal A}(F)\cong {\cal A}(G)$ if and only if $F=G$. 
\end{enumerate}
\end{Thm}

From (\ref{eqn:ostar}) and Theorem \ref{Thm:mainfive},
we obtain the following one-to-one correspondences:
\[
\begin{array}{llc}
{\cal P}({\bf N})\quad &\Leftrightarrow\quad &
\mbox{the set of all closed two-sided ideals of the C$^{*}$-algebra }\co{*},\\
\cup &&\cup\\
{\cal P}({\bf P})\quad &\Leftrightarrow\quad &
\mbox{the set of all closed biideals of the C$^{*}$-bialgebra }\co{*}.\\
\end{array}
\]
Especially,
the following anti-isomorphism between two ordered sets holds:
%
% Equation 1.24
%
\begin{equation}
\label{eqn:orderthree}
{\cal P}({\bf P})\ni F\longmapsto {\cal I}(F)\in {\rm BI}(\co{*})
\end{equation}
where ${\rm BI}(\co{*})$ is as in Corollary \ref{cor:order}.
%
% Remark 1.9
%
\begin{rem}
\label{rem:conclusion}
{\rm 
Theorem \ref{Thm:mainfive} is derived from
not only the direct sum in (\ref{eqn:cuntbi}) but also
definitions of the special comultiplication $\delp$ in (\ref{eqn:dpone})
and $\varphi_{n,m}$ in (\ref{eqn:embeddingone}).
From (\ref{eqn:dptwo}),
we see that the operation $\delp$ 
corresponds to the factorization of a natural number to two factors
in all instances.
}
\end{rem}
As a relation between lattices and operator algebras,
a lattice of von Neumann algebras is considered in \cite{Araki}.

In $\S$ \ref{section:second},
we will show a general theory of construction of C$^{*}$-bialgebra
from a family of C$^{*}$-algebras and $*$-homomorphisms.
In $\S$ \ref{subsection:secondthree}, main theorems are proved.

%%%%%%%%%%%%%%%%%%%%%%%%%%%%%%%%%%%%%%%%%%%%%%%%%%%%%%%%%%%%%%%%%%%%%%%%%%%%%%%%%%%
%
% Section 2
%
\sftt{Proofs of main theorems}
\label{section:second}
In this section, we prove main theorems in $\S$ \ref{section:first}.
%%%%%%%%%%%%%%%%%%%%%%%%%%%%%%%%%%%%%%%%%%%%%%%%%%%%%%%%%%%%%%%%%%
%
% subsection 2.1
%
\ssft{C$^{*}$-weakly coassociative system}
\label{subsection:secondone}
In order to prove main theorems in $\S$ \ref{section:first},
we show general statements about C$^{*}$-subialgebras and biideals.
We review C$^{*}$-weakly coassociative system in $\S$ 3 of \cite{TS02}.
%
% Definition 2.1
% 
\begin{defi}
\label{defi:axiom}
Let $\sem$ be a monoid with the unit $e$.
A data $(\{A_{a}:a\in \sem\},\{\varphi_{a,b}:a,b\in \sem\})$
is a C$^{*}$-weakly coassociative system (= C$^{*}$-WCS) over $\sem$ if 
$A_{a}$ is a unital C$^{*}$-algebra for $a\in \sem$
and $\varphi_{a,b}$ is a unital $*$-homomorphism
from $A_{ab}$ to $A_{a}\otimes A_{b}$
for $a,b\in \sem$ such that
\begin{enumerate}
%(i)
\item
for all $a,b,c\in \sem$, the following holds:
%
% Equation 2.1
%
\begin{equation}
\label{eqn:wcs}
(id_{a}\otimes \varphi_{b,c})\circ \varphi_{a,bc}
=(\varphi_{a,b}\otimes id_{c})\circ \varphi_{ab,c}
\end{equation}
where $id_{x}$ denotes the identity map on $A_{x}$ for $x=a,c$,
%(ii)
\item
there exists a counit $\vep_{e}$ of $A_{e}$ 
such that $(A_{e},\varphi_{e,e},\vep_{e})$ is a counital C$^{*}$-bialgebra,
%(iii)
\item
$\varphi_{e,a}(x)=I_{e}\otimes x$ and
$\varphi_{a,e}(x)=x\otimes I_{e}$ for $x\in A_{a}$ and $a\in \sem$.
\end{enumerate}
\end{defi}

\noindent
The system $(\{\con:n\in {\bf N}\},\{\varphi_{n,m}:n,m\in {\bf N}\})$
in (\ref{eqn:embeddingone}) is a C$^{*}$-WCS.
As for the other example of C$^{*}$-WCS, see $\S$ 1.3 of \cite{TS05}.

%
% Theorem 2.2
% 
\begin{Thm}(\cite{TS02}, Theorem 3.1)
\label{Thm:subthree}
Let $(\{A_{a}:a\in \sem\},\{\varphi_{a,b}:a,b\in \sem\})$ be a C$^{*}$-WCS 
over a monoid $\sem$.
Assume that $\sem$ satisfies that 
%
% Equation 2.2
%
\begin{equation}
\label{eqn:finiteness}
\#{\cal N}_{a}<\infty \mbox{ for each }a\in \sem
\end{equation}
where ${\cal N}_{a}\equiv\{(b,c)\in \sem\times \sem:\,bc=a\}$.
Then there exists a comultiplication $\delp$ and a counit $\vep$ of
the C$^{*}$-algebra 
%
% Equation 2.3
%
\begin{equation}
\label{eqn:astar}
A_{*}\equiv  \oplus \{A_{a}:a\in \sem\}
\end{equation}
such that $(A_{*},\delp,\vep)$ is a strictly proper counital C$^{*}$-bialgebra.
\end{Thm}

\noindent
We call $(A_{*},\Delta_{\varphi},\vep)$ in 
Theorem \ref{Thm:subthree} by a (counital)
{\it C$^{*}$-bialgebra} associated with 
$(\{A_{a}:a\in \sem\},\{\varphi_{a,b}:a,b\in \sem\})$.

%%%%%%%%%%%%%%%%%%%%%%%%%%%%%%%%%%%%%%%%%%%%%%%%%%%%%%%%%%
%
% subsection 2.2
%
\ssft{C$^{*}$-bialgebra associated with submonoid}
\label{subsection:secondtwo}
In this subsection, we show general results
for C$^{*}$-subbialgebras and biideals of the C$^{*}$-bialgebra $A_{*}$
in $\S$ \ref{subsection:secondone}.
Let $(\{A_{a}:a\in \sem\},\{\varphi_{a,b}:a,b\in \sem\})$ 
be a C$^{*}$-WCS over a monoid $\sem$.
For a subset $S$ of $\sem$,
define
%
% Equation 2.4
%
\begin{equation}
\label{eqn:astars}
A_{*}(S)\equiv \oplus\{A_{a}:a\in S \}\quad(S\ne \emptyset),\quad
A_{*}(\emptyset)\equiv 0.
\end{equation}
Then $A_{*}(S)$ is a closed two-sided ideal of the C$^{*}$-algebra $A_{*}$
in (\ref{eqn:astar}).

We consider the condition (\ref{eqn:biideal}) with respect to $A_{*}$.
If $H$ is a submonoid $\sem$,
then $(\{A_{a}:a\in H\},\{\varphi_{a,b}:a,b\in H\})$ 
is also a C$^{*}$-WCS over $H$.
If $\sem$ satisfies (\ref{eqn:finiteness}),
then we can define
%
% Equation 2.5
%
\begin{equation}
\label{eqn:astarh}
A_{*}(H)=\oplus \{A_{a}:a\in H\},\quad
\Delta_{\varphi^{(H)}}\equiv 
\oplus\{\Delta_{\varphi^{(H)}}^{(a)}:a\in H\},
\end{equation}
%
% Equation 2.6
%
\begin{equation}
\label{eqn:coala}
\Delta_{\varphi^{(H)}}^{(a)}(x)
\equiv \sum_{(b,c)\in {\cal N}_{a}(H)}\varphi_{b,c}(x)\quad(x\in A_{a},\,a\in H)
\end{equation}
where ${\cal N}_{a}(H)\equiv \{(b,c)\in H^{2}:bc=a\}$.
Then 
$A_{*}(H)$ is a closed two-sided ideal of the C$^{*}$-algebra $A_{*}$ and
$(A_{*}(H),\Delta_{\varphi^{(H)}},\vep|_{A_{*}(H)})$ is a counital C$^{*}$-bialgebra 
but $\Delta_{\varphi^{(H)}}\ne \delp|_{A_{*}(H)}$ in general.
Recall notions about monoid before Theorem \ref{Thm:mainfours}.
%
% Lemma 2.3
%
\begin{lem}
\label{lem:factoriallemma}
\begin{enumerate}
%(i)
\item
For a submonoid $H$ of $\sem$,
$A_{*}(H)$ is a C$^{*}$-subbialgebra of $A_{*}$
%$\Delta_{\varphi^{(H)}}= \delp|_{A_{*}(H)}$.
if and only if $H$ is factorial or $H=\sem$.
%(ii)
\item
If $H$ is a prime ideal of $\sem$,
then
$A_{*}(H)$ is a closed biideal of $(A_{*},\delp,\vep)$.
\end{enumerate}
\end{lem}
%
% Proof
%
\pr
(i)
We see that $A_{*}(H)=A_{*}$ if and only if $H=\sem$.
If $H$ is factorial, then 
$A_{*}(H)$ is a C$^{*}$-subbialgebra of $A_{*}$
from (\ref{eqn:coala}).
Assume that  $A_{*}(H)$ is a C$^{*}$-subbialgebra of $A_{*}$
but $A_{*}(H)\ne A_{*}$.
From this, $H\ne \sem$.
For $n\in H$, let $I_{n}$ denote the unit of $A_{n}$.
Then 
%
% Equation 2.7
%
\begin{equation}
\label{eqn:units}
\delp(I_{n})=\delp^{(n)}(I_{n})
=\sum_{(m,l)\in {\cal N}_{n}}\varphi_{m,l}(I_{n})
=\sum_{(m,l)\in {\cal N}_{n}}I_{m}\otimes I_{l}.
\end{equation}
By assumption, 
%
% Equation 2.8
%
\begin{equation}
\label{eqn:unitst}
\sum_{(m,l)\in {\cal N}_{n}}I_{m}\otimes I_{l}\in 
\delp(A_{*}(H))\subset A_{*}(H)\otimes A_{*}(H)
=\bigoplus_{n^{'},n^{''}\in H}A_{n^{'}}\otimes A_{n^{''}}.
\end{equation}
By multiplying $I_{m}\otimes I_{l}$ at both sides of (\ref{eqn:unitst}),
$I_{m}\otimes I_{l}\in 
\bigoplus_{n^{'},n^{''}\in H}I_{m}A_{n^{'}}\otimes I_{l}A_{n^{''}}$.
Hence both $m$ and $l$ must belong to $H$
because $I_{a}A_{b}=\{0\}$ when $a\ne b$.
Therefore $H$ is factorial.

\noindent
(ii)
It is sufficient to show that $A_{*}(H)$ satisfies
(\ref{eqn:biideal}) with respect to  $(A_{*},\delp,\vep)$.
For $n\in H$ and $x\in A_{n}$,
%
% Equation 2.9
%
\begin{equation}
\label{eqn:delx}
\delp(x)=\delp^{(n)}(x)=
\sum_{(m,l)\in {\cal N}_{n}}\varphi_{m,l}(x).
\end{equation}
If $(m,l)\in {\cal N}_{n}$,
then $m\in H$ or $l\in H$.
Therefore
$\varphi_{m,l}(x)\in A_{*}(H)\otimes A_{*}+A_{*}\otimes A_{*}(H)$.
This implies that $\delp(x)\in A_{*}(H)\otimes A_{*}+A_{*}\otimes A_{*}(H)$.
Since $H$ is a prime ideal of $\sem$,
the unit $e$ of $\sem$ does not belong to $H$.
Therefore $\vep(A_{*}(H))=0$.
Hence the statement holds.
\qedh

By definition, the following holds.
%
% Fact 2.4
%
\begin{fact}
\label{fact:two}
For a subset $H$ of $\sem$,
let $H^{c}$ denote $\{x\in \sem:x\not\in H\}$.
Then the following holds:
\begin{enumerate}
%(i)
\item
$H$ is a proper subsemigroup if and only if $H^{c}$ is a prime subset.
%(ii)
\item
$H$ is a proper ideal if and only if $H^{c}$ is factorial.
\end{enumerate}
\end{fact}
%
% Proof
%
\pr
(i)
Assume that $H$ is a proper subsemigroup.
Then $H^{c}\ne \emptyset$, and 
if $a,b\in\sem$ and  $a,b\in H$, then $ab\in H$.
Therefore $a\not\in H$ or $b\not\in H$ when $ab\not \in H$.
This implies that $a\in H^{c}$ or $b\in H^{c}$ when $ab\in H^{c}$.
Therefore $H^{c}$ is prime.

Assume that $H^{c}$ is prime.
Then $H\ne \emptyset$.
For $a,b\in\sem$, if 
$ab\not \in H$,
then $ab\in H^{c}$.
Since $H^{c}$ is prime, $a\in H^{c}$ or $b\in H^{c}$.
Hence $a\not \in H$ or $b\not \in H$.
From contraposition,
if $a,b\in H$, then $ab\in H$.
Since $H^{c}\ne \emptyset$,
Therefore $H$ is a proper subsemigroup.

\noindent
(ii)
Assume that $H$ is a proper ideal.
Then $H^{c}\ne \emptyset$.
For $a,b\in\sem$,
if $ab\in H^{c}$, then $ab\not\in H$.
Since $H$ is an ideal,
$a,b\not\in H$.
Hence $a,b\in H^{c}$.
Therefore $H^{c}$ is factorial.

Assume that $H^{c}$ is factorial.
Then $H\ne \emptyset$.
For $a,b\in\sem$,
if $a\in H$, then $a\not\in H^{c}$.
In addition,
if $ab\in H^{c}$, then $a$ must belong to $H^{c}$
because $H^{c}$ is factorial.
Hence $ab\not\in H^{c}$.
This implies that $ab\in H$.
In the same way, 
if $b\in H$, then $ab\in H$.
Hence $H$ is an ideal.
Since $H^{c}\ne\emptyset$, $H$ is proper.
\qedh
%
% Proposition 2.5
%
\begin{prop}
\label{prop:general}
Let $(\{A_{a}:a\in \sem\},\{\varphi_{a,b}:a,b\in \sem\})$ 
be as in Theorem \ref{Thm:subthree}.
In addition,
assume that $A_{a}$ is simple for each $a\in \sem$.
\begin{enumerate}
%(i)
\item
For a  subspace ${\cal I}$ of $A_{*}$,
${\cal I}$ is a closed biideal 
of $(A_{*},\delp,\vep)$ if and only if
\begin{enumerate}
%(a)
\item
${\cal I}=\{0\}$, or 
%(b)
\item
there exists a prime ideal ${\goth a}$ of $\sem$
such that ${\cal I}=A_{*}({\goth a})$.
\end{enumerate}
%(ii)
\item
For any closed biideal ${\cal I}$ of $(A_{*},\delp,\vep)$,
there exists a C$^{*}$-subbialgebra ${\cal B}_{{\cal I}}$ of $(A_{*},\delp,\vep)$
such that the quotient C$^{*}$-bialgebra 
$A_{*}/{\cal I}$ is isomorphic to ${\cal B}_{{\cal I}}$
and the following decomposition holds:
%
% Equation 2.10
%
\begin{equation}
\label{eqn:decompositionb}
A_{*}={\cal B}_{{\cal I}}\oplus {\cal I}
\end{equation}
where $\oplus$ means the direct sum  of two C$^{*}$-subalgebras.
\end{enumerate}
\end{prop}
%
% Proof
%
\pr
(i)
We see that $\{0\}$ is a closed biideal of $A_{*}$.
Assume that ${\cal I}$ is a nonzero closed biideal 
of $(A_{*},\delp,\vep)$.
Since ${\cal I}$ is a nonzero closed two-sided ideal of $A_{*}$ and 
$A_{a}$ is simple for each $a\in \sem$, 
we see that there exists a nonempty subset $S$ of $\sem$
such that ${\cal I}=A_{*}(S)$.
Let $n\in S$.
From the assumption (\ref{eqn:biideal}) and (\ref{eqn:units}),
%
% Equation 2.11
%
\begin{equation}
\label{eqn:summl}
\sum_{(m,l)\in {\cal N}_{n}}I_{m}\otimes I_{l}
=\delp(I_{n})\in {\cal I}\otimes A_{*}+A_{*}\otimes {\cal I}
=\bigoplus_{n^{'}\in S}A_{n^{'}}\otimes A_{*}
+A_{*}\otimes \bigoplus_{n^{''}\in S}A_{n^{''}}.
\end{equation}
By multiplying $I_{m}\otimes I_{l}$ at both sides of (\ref{eqn:summl}),
%
% Equation 2.12
%
\begin{equation}
\label{eqn:times}
I_{m}\otimes I_{l}
\in\bigoplus_{n^{'}\in S}I_{m}A_{n^{'}}\otimes A_{l}
+A_{m}\otimes \bigoplus_{n^{''}\in S}I_{l}A_{n^{''}}.
\end{equation}
If $m\not\in S$ in (\ref{eqn:times}), then 
$\bigoplus_{n^{'}\in S}I_{m}A_{n^{'}}\otimes A_{l}=\{0\}$
and $I_{m}\otimes I_{l}
\in A_{m}\otimes \bigoplus_{n^{''}\in S}I_{l}A_{n^{''}}$.
Hence $l$ must belong to $S$.
By the same token,
if $l\not\in S$ in (\ref{eqn:times}), then $m$ belongs to $S$.
This implies that 
if $n=ml$, then $m\in S$ or $l\in S$.
Therefore $S$ is prime. 
If $S=\sem$, then $S$ is an ideal of $\sem$.
If $S\ne \sem$, then 
$S^{c}$ is a subsemigroup of $\sem$
from Fact \ref{fact:two}(i).
From the assumption (\ref{eqn:biideal}),
$\vep(A_{*}(S))=0$.
Hence the unit $e$ does not belong to $S$.
Therefore $S^{c}$ is a submonoid of $\sem$.
Since ${\cal I}$ is a closed biideal of $A_{*}$,
$A_{*}/{\cal I}=A_{*}/A_{*}(S)$ is isomorphic to
$A_{*}(S^{c})$ as a C$^{*}$-bialgebra
such that the natural projection
from $A_{*}$ to $A_{*}/{\cal I}$ is a C$^{*}$-bialgebra morphism. 
This implies that 
$A_{*}(S^{c})$ is a  C$^{*}$-subbialgebra of $A_{*}$.
From Lemma \ref{lem:factoriallemma}(i), $S^{c}$ is factorial.	
Therefore $S^{c}$ is a factorial submonoid of $\sem$.
From Fact \ref{fact:inverse}, $S$ is a prime ideal of $\sem$.

The inverse statement holds 
from Lemma \ref{lem:factoriallemma}(ii).

\noindent
(ii)
In the proof of (i),
let ${\cal B}_{{\cal I}}\equiv A_{*}(S^{c})$.
Then the statement holds.
\qedh

From Proposition \ref{prop:general}(i) and (\ref{eqn:astars}), 
the following holds.
%
% Corollary 2.6
%
\begin{cor}
\label{cor:ordertwo}
Let $A_{*}$ and $\sem$ be as in Proposition \ref{prop:general}.
Let ${\rm Spec}\sem$ denote the set of all prime ideals of the monoid $\sem$.
Define $\overline{{\rm Spec}\sem}\equiv {\rm Spec}\sem\cup \{\emptyset\}$
and let ${\rm BI}(A_{*})$ denote the set of all closed biideals of $A_{*}$.
Then the following isomorphism between two ordered sets holds:
%
% Equation 2.13
%
\begin{equation}
\label{eqn:ordertwo}
\overline{{\rm Spec}\sem}\ni x\longmapsto A_{*}(x)\in {\rm BI}(A_{*})
\end{equation}
where the order of both 
$\overline{{\rm Spec}\sem}$ and ${\rm BI}(A_{*})$
is taken as the inclusion of subsets.
\end{cor}
%
% Proof
%
\pr
If $x=\emptyset$,
then $A_{*}(\emptyset)=\{0\}$ by (\ref{eqn:astars}).
If $x\ne \emptyset$,
then 
the map $x\mapsto A_{*}(x)$ is bijective
from ${\rm Spec}\sem$ to ${\rm BI}(A_{*})\setminus \{\{0\}\}$
from Proposition \ref{prop:general}(i)-(a).
The order structure follows from (\ref{eqn:astars}).
\qedh

Remark that we {\it do not} assume that $\sem$ is abelian
in this subsection.

%%%%%%%%%%%%%%%%%%%%%%%%%%%%%%%%%%%%%%%%%%%%
%
% Subsection 2.3
%
\ssft{Proofs of main theorems}
\label{subsection:secondthree}
In this subsection, we prove main theorems.
\\
\\
{\it Proof of Theorem \ref{Thm:maintwo}.}
Applying Theorem \ref{Thm:subthree}
for ${\sf M}=H$ and $A_{n}=\con$ with $n\in H$, the statement holds.
\qedh
\\
{\it Proof of Theorem \ref{Thm:mainfours}.}
Applying Proposition \ref{prop:general}
for ${\sf M}={\bf N}$ and $A_{n}=\con$ with $n\in {\bf N}$, the statement holds.
\qedh
\\
{\it Proof of Corollary \ref{cor:order}.}
Applying Corollary \ref{cor:ordertwo}
for ${\sf M}={\bf N}$ and $A_{n}=\con$ with $n\in {\bf N}$, 
the statement holds. 
\qedh
\\
{\it Proof of Fact \ref{fact:inverse}.}
From Fact \ref{fact:two},
the statement holds. 
\qedh
\\
{\it Proof of Theorem \ref{Thm:mainfive}.}
(i)
The proof is already given in statements before 
Theorem \ref{Thm:mainfive}.
We summarize them again as follows.
Let ${\cal J}$ be a nonzero closed biideal of $\co{*}$.
From Theorem \ref{Thm:mainfours}(i),
this is equivalent that 
there exists a prime ideal ${\goth a}$ of ${\bf N}$
such that ${\cal J}=\co{*}({\goth a})$.
From Fact \ref{fact:inverse},
this is equivalent that 
there exists a factorial submonoid $H$ of ${\bf N}$
such that ${\cal J}=\co{*}(H^{c})$.
From (\ref{eqn:powers}),
this is equivalent that 
there exists $F\in {\cal P}({\bf P})$
such that ${\cal J}=\co{*}([F]^{c})={\cal I}(F)$.

\noindent
(ii)
It is sufficient to show
${\cal P}({\bf P})\ni F\mapsto {\cal A}(F)$
is a lattice isomorphism.
For this purpose,
we show the following statements:

\def\labelenumi{\theenumi}
\def\theenumi{(\alph{enumi})}

\begin{enumerate}
%(i)
\item
${\cal A}(\emptyset)=\co{1}$,
${\cal A}({\bf P})=\co{*}$.
%(ii)
\item
For $F,G\in {\cal P}({\bf P})$,
${\cal A}(F)\cap {\cal A}(G)={\cal A}(F\cap G)$ and
${\cal A}(F)\vee {\cal A}(G)={\cal A}(F\cup G)$.
%(iii)
\item
For $F,G\in {\cal P}({\bf P})$,
${\cal A}(F)\cong {\cal A}(G)$ if and only if $F=G$.
\end{enumerate}

\renewcommand{\thesection}{\arabic{section}}

From  $[{\bf P}]={\bf N}$, $[\emptyset]=\{1\}$ and (\ref{eqn:af}),
the statement (a) follows.
For $F,G\in {\cal P}({\bf P})$,
%
% Equation 2.13
%
\begin{equation}
\label{eqn:cap}
{\cal A}(F)\cap {\cal A}(G)=\co{*}([F])\cap \co{*}([G])
=\co{*}([F]\cap [G])=\co{*}([F\cap G])={\cal A}(F\cap G).
\end{equation}
By the same token,
${\cal A}(F)\vee {\cal A}(G)={\cal A}(F\cup G)$.
By using these and (i), other statements of (b) follow.
For $F,G\in {\cal P}({\bf P})$,
${\cal A}(F)\cong {\cal A}(G)$ if and only if 
$\co{*}([F])\cong\co{*}([G])$.
This is equivalent to $[F]=[G]$ from the statement after (\ref{eqn:ostar}).
Since ${\bf N}$ is the free abelian semigroup generated by ${\bf P}$,
this is equivalent to $F=G$.
Hence (c) holds.
\qedh

%\ww
%{\bf Acknowledgement:}
%The author would like to express his sincere thanks to Izumi Ojima 
%for his interest in this topic and raising the above question.

%%%%%%%%%%%%%%%%%%%%%%%%%%%%%%%%%%%%%%%%%%%%%%%%%%%%%%%%%%%%%%%%%%%%%
%
% Reference 
%

%
\end{document}